\documentclass[12pt]{article}
\usepackage{mathrsfs}
\usepackage{amsfonts}
\usepackage{amsthm,amsmath,amssymb,anysize,longtable}

\newtheorem{theorem}{Theorem}
\newtheorem{lemma}[theorem]{Lemma}

\newtheorem{definition}[theorem]{Definition}
\newtheorem{corollary}[theorem]{Corollary}
\usepackage{xcolor}

\usepackage[numbers,sort&compress]{natbib}
\usepackage{amsmath}
\allowdisplaybreaks[3]
\marginsize{4.2cm}{4.2cm}{3.6cm}{3cm}
\numberwithin{equation}{section}

\author{\textsc{Xuelian Guo and Liming Tang\footnote{corresponding author}}\\
\small{School of Mathematical Sciences}\\
\small{Harbin Normal University}\\
\small{150025 Harbin, China}\\
\small{E-mail: limingtang@hrbnu.edu.cn}}

\date{ }
\date{ }
\usepackage[symbol]{footmisc}

\begin{document}

\thispagestyle{empty}

\noindent{\Large
Derivations, $2$-local derivations, biderivations and automorphisms of generalized loop Heisenberg-Virasoro algebra}
\footnote{
The work is supported by the NSF of Hei Longjiang Province (No. LH2024A014).
}

	\bigskip
	
	 \bigskip

\begin{center}	
	{\bf
		
    Qingyan Ren\footnote{School of Mathematical Sciences, Harbin Normal University, 150025 Harbin, China;\ renqingyan@stu.hrbnu.edu.cn},
Liming Tang\footnote{School of Mathematical Sciences, Harbin Normal University, 150025 Harbin, China; \ limingtang@hrbnu.edu.cn}\footnote{corresponding author}}
\end{center}

 \begin{quotation}
{\small\noindent \textbf{Abstract}:
In this paper, the generalized loop Heisenberg-Virasoro algebra is introduced. Firstly, we determine the derivations on the generalized loop Heisenberg-Virasoro algebra. Then we show that all $2$-local derivations are derivations for the generalized loop Heisenberg-Virasoro algebra. Furthermore, we prove the biderivations on the generalized loop Heisenberg-Virasoro algebra are inner biderivations and give their applications. Finally, the automorphism groups on the generalized loop Heisenberg-Virasoro algebra are presented.

\medskip
 \vspace{0.05cm} \noindent{\textbf{Keywords}}:
 generalized loop Heisenberg-Virasoro algebras; derivations; $2$-local derivations; biderivations; automorphism groups

\medskip

\vspace{0.05cm} \noindent \textbf{Mathematics Subject Classification
2020}: 17B40, 17B65, 17B68}
\end{quotation}
 \medskip

\section*{Introduction}

The Heisenberg-Virasoro algebra is an important infinite dimensional Lie algebra and it plays a vital role in numerous fields of mathematics and physics. The Heisenberg-Virasoro algebra is the universal central extension of the Lie algebra of differential operators of degree at most 1 defined on the circle, and has been extensively studied in mathematical and physical literature.

\smallskip

As is well known, structure theory is essential in the research of Lie algebras. Likewise, the authors in \cite{ref5,ref6,ref8,ref9,ref18} determined the derivation algebra and automorphism group on the generalized loop Virasoro algebra, the generalized loop Schrödinger–Virasoro algebras and the generalized Heisenberg-Virasoro algebra, etc. The authors in \cite{ref4} obtain $2$-local derivations on the generalized Witt algebras. The authors in \cite{ref19} proved that all 2-local derivations on the Witt algebra as well as on the positive Witt algebra are derivations, and provided an example of infinite-dimensional Lie algebra with a 2-local derivation which is not a derivation. The authors in \cite{ref2} determined biderivations on the three Lie algebras of maximal class with one-dimensional homogeneous components, and characterized their local and 2-local derivations. The authors in \cite{ref10,ref11,ref12,ref17} gave the biderivations on the twisted Heisenberg–Virasoro algebra, the Schrödinger-Virasoro algebra and the W-algebras. The authors in \cite{ref20} described all biderivations of a semisimple Lie algebra. The authors in \cite{ref1} described the structure of the irreducible highest weight modules for the twisted Heisenberg-Virasoro Lie algebra at level zero. The authors in \cite{ref14} constructed two classes of non-weight modules over the twisted Heisenberg-Virasoro algebra and presented the necessary and sufficient conditions under which modules in these two classes are irreducible and isomorphic. The authors in \cite{ref13}  gave sufficient and necessary conditions for these non-weight modules on the mirror Heisenberg-Virasoro algebra to be irreducibie. The authors in \cite{ref15} studied irreducible modules over the mirror Heisenberg-Virasoro algebra and determined irreducible modules with finite-dimensional weight spaces are either an irreducible highest or an irreducible lowest weight module, or an irreducible module of the intermediate series. The authors in \cite{ref16} determined the necessary and sufficient conditions for the tensor products of irreducible highest weight modules and irreducible modules of intermediate series over the mirror Heisenberg-Virasoro algebra to be irreducible by using “shifting technique”.

\smallskip

In this paper, inspired by \cite{ref7}, we introduce the generalized loop Heisenberg-Virasoro algebra. 
The aim is to determine the derivations, $2$-local derivations, biderivations and automorphism groups of the generalized loop Heisenberg-Virasoro algebra by use of the analogous methods with those of respective Lie algebras in \cite{ref6,ref8,ref11,ref19,ref9}. The structure of this paper is as follows: Firstly, based on its relevant properties, the derivations on the generalized loop Heisenberg-Virasoro algebra are determined. On this basis, it is proved that all $2$-local derivations of this algebra are derivations. Secondly, based on the fact that it is a perfect Lie algebra, we conclude that all biderivations are inner biderivations. Then we prove the necessary and sufficient conditions for its commutating linear maps and any commutative $post$-Lie algebra structure on this algebra is trivial. Finally, we define five types of automorphism groups of the algebra and prove that any automorphism group on it can be represented the direct product of the five types of automorphisms.

\section{Preliminaries}
Throughout this paper, the ground field $\mathbb{F}$ is an algebraically closed field of characteristic zero and all vector spaces, algebras are over $\mathbb{F}$. $\Gamma$ represents an abelian additive group, and $\mathbb{Z}$ represents the set of integers. Let $\mathrm{span}\{X\}$ denote the vector space generated by the vectors in $\{X\}$.

\begin{definition}
Generalized loop Heisenberg-Virasoro algebra is a Lie algebra with basis $\left\{L_{\alpha,i},H_{\beta,j}\,|\,\alpha,\beta\in\Gamma,\,i,j\in\mathbb{Z}\right\}$, subject to the following nontrivial Lie brackets:
\begin{longtable}{lcl}
$[L_{\alpha,i},L_{\beta,j}]$&$=$&$(\alpha-\beta)L_{\alpha+\beta,i+j},$ \\
$[L_{\alpha,i},H_{\beta,j}]$&$=$&$-\beta H_{\alpha+\beta,i+j},$ \\
$[H_{\alpha,i},H_{\beta,j}]$&$=$&$0.$
\end{longtable}
\end{definition}

In the sequel, denoted by $\mathcal{L}(\Gamma)$ generalized loop Heisenberg-Virasoro algebra. It is worth noting that loop Heisenberg-Virasoro algebra in \cite{ref7} is exactly $\mathcal{L}(\Gamma)$ when $\mathbb{Z}=\Gamma$. Let  $\mathcal{LV}$ be subalgebra of $\mathcal{L}(\Gamma)$ spanned by $\{L_{\alpha,i}\,|\,\alpha\in\Gamma,\,i\in\mathbb{Z}\}$. Then $\mathcal{LV}$ is actually isomorphic to the centerless generalized loop Virasoro algebra (see \cite{ref6}). In addition, $\mathcal{L}(\Gamma)$ contains the generalized Heisenberg-Virasoro algebra (see \cite{ref9})
$$\mathcal{HV}=\mathrm{span}_\mathbb{F}\{L_{\alpha,0},H_{\beta,0}\,|\,\alpha,\beta\in\Gamma\}.$$

Let $L$ be a Lie algebra and $V$ be $L$-module. A linear map $\mathcal{D}:L\rightarrow V$ is called a derivation of $\mathcal{L}(\Gamma)$ if
$$\mathcal{D}([x,y])=x\cdot\mathcal{D}(y)+y\cdot\mathcal{D}(x)\, \mbox{for all}\,  x,y\in L.$$
A linear map $\mathrm{ad}: x\mapsto x\cdot v$ for all $x\in L$ and all $v\in V$ is called the inner derivation. Otherwise others are called outer derivations.

Let $\mathrm{Der}(\mathcal{L}(\Gamma))$ and $\mathrm{Inn}(\mathcal{L}(\Gamma))$ be the vector spaces consisting of all derivations and all inner derivations on $\mathcal{L}(\Gamma)$, respectively. Then $$H^1(\mathcal{L}(\Gamma))\cong\mathrm{Der}(\mathcal{L}(\Gamma))/\mathrm{Inn}(\mathcal{L}(\Gamma))$$ is the first cohomology group of $\mathcal{L}(\Gamma)$ (see \cite{ref8}).

Since $L_{0,0}$ is semisimple, then the $\Gamma$-grading on $\mathcal{L}(\Gamma)$ is given, $$\mathcal{L}(\Gamma)=\bigoplus_{\mu \in \Gamma}\mathcal{L}(\Gamma)_\mu,$$
where
$$\mathcal{L}(\Gamma)_\mu=\left\{x\in \mathcal{L}(\Gamma)\,|\,[L_{0,0},x]=-\mu x\right\}=\mathrm{span}_\mathbb{F}\left\{L_{\mu,j},H_{\mu,j}\,|\,j\in\mathbb{Z}\right\}.$$
Recall that a derivation $\mathcal{D}\in\mathrm{Der}\mathcal{L}(\Gamma)$ is of degree $\gamma\in\Gamma$ if $D(\mathcal{L}(\Gamma)_\alpha)\subset\mathcal{L}_{\alpha+\gamma}$ for all $\alpha\in\Gamma$. Let $(\mathrm{Der}\mathcal{L}(\Gamma))_\gamma$ be the space of all derivations of degree $\gamma$.

\section{Derivations on $\mathcal{L}(\Gamma)$}
\begin{lemma}
$\mathrm{Der}\mathcal{L}(\Gamma)=(\mathrm{Der}\mathcal{L}(\Gamma))_0+\mathrm{ad}\mathcal{L}(\Gamma)$.

\begin{proof}
(\romannumeral1) First we show that $\mathrm{Der}\mathcal{L}(\Gamma)=\sum_{\gamma\in\Gamma}(\mathrm{Der}\mathcal{L}(\Gamma))_{\gamma}.$

For all $\mathcal{D}\in\mathrm{Der}\mathcal{L}(\Gamma)$ and $x=\sum_{\alpha\in\Gamma}x_\alpha$, where $x_\alpha\in\mathcal{L}(\Gamma)_\alpha$, assume that $\mathcal{D}(x_\alpha)=\sum_{\beta\in\Gamma}y_\beta$. We define $\mathcal{D}_\gamma(x_\alpha)=y_{\alpha+\gamma}$. It is easy to prove that
$$\mathcal{D}_\gamma([x_\alpha,x_\beta])=[\mathcal{D}_\gamma(x_\alpha),x_\beta]+[x_\alpha,\mathcal{D}_\gamma(x_\beta)],\, \mbox{for all}\, x_\alpha\in\mathcal{L}(\Gamma)_\alpha \mbox{and all}\,x_\beta\in\mathcal{L}(\Gamma)_\beta.$$
Then $\mathcal{D}_\gamma$ is a derivation on $\mathcal{L}(\Gamma)$. We can conclude that
$$\mathcal{D}(x)=\sum_{\alpha\in\Gamma}\mathcal{D}(x_\alpha)=\sum_{\alpha\in\Gamma}\sum_{\beta\in\Gamma}y_\beta=\sum_{\gamma\in\Gamma}\mathcal{D}_\gamma(x).$$
Then $\mathcal{D}=\sum_{\gamma\in\Gamma}\mathcal{D}_{\gamma}$. That is for any derivation $\mathcal{D}\in\mathrm{Der}\mathcal{L}(\Gamma)$, one has $\mathcal{D}=\sum_{\gamma\in\Gamma}\mathcal{D}_{\gamma}$, where $\mathcal{D}_\gamma\in(\mathrm{Der}\mathcal{L}(\Gamma))_\gamma$,
such that for all $x\in\mathcal{L}(\Gamma)$, only finitely $\mathcal{D}_{\gamma}(x)\neq0$.

(\romannumeral2) Then we show  that for all $\gamma\neq0 $, $\mathcal{D}\in(\mathrm{Der}\mathcal{L}(\Gamma))_\gamma$ is an inner derivation.

For all $x\in\mathcal{L}(\Gamma)_\mu$, on the one hand, it has
$$\mathcal{D}([L_{0,0},x])=[\mathcal{D}(L_{0,0}),x]{-(\mu+\gamma)}\mathcal{D}(x),$$
on the other hand,
\begin{longtable}{lcl}
$\mathcal{D}([L_{0,0},x])$&$=$&$\mathcal{D}(-\mu x)$\\
$ $&$=$&$-\mu\mathcal{D}(x)$,
\end{longtable}
then
\begin{longtable}{lcl}
$\mathcal{D}(x)$&$=$&$[\frac{1}{\gamma}\mathcal{D}(L_{0,0}),x]$\\
$ $&$=$&${\rm{ad}}_{\gamma^{-1}\mathcal{D}(L_{0,0})}(x)$.
\end{longtable}
Thus for all $\gamma\neq 0$, $\mathcal{D}=\mathrm{ad}_{\gamma^{-1}\mathcal{D}(L_{0,0})}$ is an inner derivation.

\end{proof}
\end{lemma}

By Lemma 2, if we  determine $\mathrm{Der}\mathcal{L}(\Gamma)$, only needs to determine $(\mathrm{Der}\mathcal{L}(\Gamma))_0$. For a map $f$ from a set $A$ to another set $B$, let $f_a$ denote the image of $a\in A$ under $f$. For all $$\phi\in \mathrm{Hom}_\mathbb{Z}(\Gamma,\mathbb{F}[t,t^{-1}]), g\in\mathbf{g}(\Gamma), b\in\mathbb{F}[t,t^{-1}], \rho\in\mathbb{F}[t,t^{-1}]\frac{d}{dt}$$ and $$\mathbf{g}(\Gamma)=\left\{g:\Gamma\rightarrow\mathbb{F}[t,t^{-1}]\,|\,(\alpha-\beta)g_{\alpha+\beta}=\alpha g_\alpha-\beta g_\beta,\,\alpha,\beta\in\Gamma\right\},$$ where $\mathbb{F}[t,t^{-1}]\frac{d}{dt}$ is a derivation algebra of $\mathbb{F}[t,t^{-1}]$.  The following linear maps $D_\phi,\,D_g,\,D_b,\,D^\rho$ on $\mathcal{L}(\Gamma)$ are given by, respectively:
\begin{align*}
D_\phi(L_{\alpha,i})&=\phi(\alpha)L_{\alpha,i}, &D_\phi(H_{\alpha,i})&=\phi(\alpha)H_{\alpha,i},\\
D_g(L_{\alpha,i})&=g_\alpha H_{\alpha,i}, &D_g(H_{\alpha,i})&=0,\\
D_b(L_{\alpha,i})&=0, &D_b(H_{\alpha,i})&=bH_{\alpha,i},\\
D^\rho(L_{\alpha,i})&=L_\alpha\rho(t^i), &D^\rho(H_{\alpha,i})&=H_\alpha\rho(t^i),\, \mbox{for all}\,\alpha\in\Gamma \,\mbox{and all}\,i\in\mathbb{Z}.
\end{align*}

It is easy to see that these operators defined above are all homogeneous derivations of degree zero. Let
\begin{longtable}{lcl}
$ $&&$\mathcal{D}_{\mathrm{Hom}_\mathbb{Z}(\Gamma,\mathbb{F}[t,t^{-1}])}\,:=\,\left\{D_\phi\,|\,\phi\in\mathrm{Hom}_\mathbb{Z}(\Gamma,\mathbb{F}[t,t^{-1}])\right\}$,\\
$ $&&$\mathcal{D}_{\mathbf{g}(\Gamma)}\,:=\,\left\{D_g\,|\,g\in\mathbf{g}(\Gamma)\right\}$,\\
$ $&&$\mathcal{D}_{\mathbb{F}[t,t^{-1}]}\,:=\,\left\{D_b\,|\,b\in\mathbb{F}[t,t^{-1}]\right\}$,\\
$ $&&$\mathcal{D}_{\mathbb{F}[t,t^{-1}]\frac{d}{dt}}\,:=\,\{D^\rho\,|\,\rho\in\mathbb{F}[t,t^{-1}]\frac{d}{dt}\}$
\end{longtable}
be the corresponding subspaces of $\mathrm{Der}\mathcal{L}(\Gamma)$, respectively. We brief note $\mathbf{g}(\Gamma)$ by $\mathbf{g}$ in the sequel.

\begin{lemma}
Suppose $D\in(\mathrm{Der}\mathcal{L}(\Gamma))_0$, then there exist $f\in\mathrm{Hom}_\mathbb{Z}(\Gamma,\mathbb{F}[t,t^{-1}]),\,g\in\mathbf{g}$  such that
\begin{longtable}{lcl}
$D(L_{\alpha,i})=f_\alpha L_{\alpha,i}+g_\alpha H_{\alpha,i}$ for all $\alpha\in\Gamma$ and all $i\in\mathbb{Z}.$
\end{longtable}

\begin{proof}
Let $D(L_{\alpha,i})=f_{\alpha,i}L_{\alpha,i}+g_{\alpha,i}H_{\alpha,i}$ where $f_{\alpha,i},g_{\alpha,i}\in\mathbb{F}[t,t^{-1}]$.
On the one hand, we have
$$D([L_{\alpha,i},L_{\beta,j}])=(\alpha-\beta)(f_{\alpha+\beta,i+j}L_{\alpha+\beta,i+j}+g_{\alpha+\beta,i+j}H_{\alpha+\beta,i+j}).$$
On the other hand,
\begin{longtable}{lcl}
$D([L_{\alpha,i},L_{\beta,j}])$&$=$&$[D(L_{\alpha,i}),L_{\beta,j}]+[L_{\alpha,i},D(L_{\beta,j})]$\\
$
$&$=$&$(\alpha-\beta)f_{\alpha,i}L_{\alpha+\beta,i+j}+\alpha g_{\alpha,i}H_{\alpha+\beta,i+j}+(\alpha-\beta)f_{\beta,j}L_{\alpha+\beta,i+j}$\\
$ $&&$-\beta g_{\beta,j}H_{\alpha+\beta,i+j}$.
\end{longtable}
Then we obtain\begin{longtable}{lcl}
$ $&&$f_{\alpha+\beta,i+j}=f_{\alpha,i}+f_{\beta,j},\,\alpha\neq\beta$,\\
$ $&&$(\alpha-\beta)g_{\alpha+\beta,i+j}=\alpha g_{\alpha,i}-\beta g_{\beta,j}$ for all $\alpha,\beta\in\Gamma$ and all $i,j\in\mathbb{Z}$.
\end{longtable}
Let $\alpha=0,i=0$. Then $$f_{0+\beta,0+j}=f_{0,0}+f_{\beta,j}.$$ Thus $f_{0,0}=0$. Let $\beta=-\alpha\neq0,i=j=0$. Then $$0=f_{0,0}=f_{\alpha,0}+f_{-\alpha,0},\,\alpha\neq0.$$
Let $0\neq\eta\in\Gamma$ satisfying $\eta\neq\alpha,\,-\eta\neq\beta,\alpha+\eta\neq\beta-\alpha$. Then
\begin{longtable}{lcl}
$ $&&$f_{\alpha+\beta,i+j}=f_{\alpha+\eta,i}+f_{\beta-\eta,j}$,\\
$ $&&$f_{\alpha+\eta,i}=f_{\alpha,i}+f_{\eta,0}$,\\
$ $&&$f_{\beta-\eta.j}=f_{\beta,j}+f_{-\eta,0}$.
\end{longtable}
That is, $f_ {\alpha+\beta,i+j}=f_{\alpha,i}+f_{\beta,j}$ for all $\alpha,\beta\in\Gamma$ and all $i,j\in\mathbb{Z}.$
Then $f_{\alpha,i}=f_{\alpha,0}+f_{0,i}$. Let $\rho=tf_{0,1}\frac{d}{dt}$. We have $$D^\rho(L_{0,1})=L_0(tf_{0,1}\frac{d}{dt})(t)=L_0tf_{0,1}.$$
Replace $D$ by $D-D^\rho$. Let $D(L_{0,1})=0$. Then $f_{0,1}=0$ and $f_{0,i}=if_{0,1}=0$, thus $f_{\alpha,i}=f_{\alpha,0}+f_{0,i}=f_{\alpha,0}$ for all $\alpha\in\Gamma$ and all $i\in\mathbb{Z}.$
For convenience, denote $f_{\alpha,i}$ as $f_\alpha$. Then $$f_{\alpha+\beta}=f_\alpha+f_\beta,\,f\in\mathrm{Hom}_\mathbb{Z}(\Gamma,\mathbb{F}[t,t^{-1}]).$$
Let $\alpha=0,j=0$. Then we have $$g_{\beta,i}=g_{\beta,0}.$$
Let $\alpha=-\beta\neq0$. Then we have
\begin{longtable}{lcl}
$2\beta g_{0,i+j}$&$=$&$\beta g_{\beta,j}+\beta g_{-\beta,i}$\\
$ $&$=$&$\beta g_{\beta,0}+\beta g_{-\beta,0}$\\
$ $&$=$&$2\beta g_{0,0}$\\
$ $&$=$&$0$.
\end{longtable}
Thus $g_{\beta,i}=g_{\beta,0}-g_{0,i}=g_{\beta,0}$ for all $\beta\in\Gamma$ and all $i\in\mathbb{Z}.$
For convenience, denote $g_\beta$ by $g_{\beta,i}$. Then $(\alpha-\beta)g_{\alpha+\beta}=\alpha g_\alpha-\beta g_\beta$ for all $g\in\mathbf{g}(\Gamma).$
\end{proof}
\end{lemma}

\begin{lemma}
Let $D\in(\mathrm{Der}\mathcal{L}(\Gamma))_0$. Then there exists  $b\in\mathbb{F}[t,t^{-1}]$ such that $D(H_{\alpha,i})=(f_\alpha+b)H_{\alpha,i}$ for all $\alpha\in\Gamma$ and all $i\in\mathbb{Z}.$

\begin{proof}
For all $\alpha\in\Gamma$ and all $i\in\mathbb{Z}$, assume that $$D(H_ {\alpha,i})=\phi_{\alpha,i}L_{\alpha,i}+\psi_{\alpha,i}H_{\alpha,i},\,\phi_{\alpha,i},\psi_{\alpha,i}\in\mathbb{F}[t,t^{-1}].$$
On the one hand,
\begin{longtable}{lcl}
$D([L_{-\alpha,-i},H_{\alpha,i}])$&$=$&$D(-\alpha H_{0,0})$\\
$ $&$=$&$-\alpha(\phi_{0,0}L_{0,0}+\psi_{0,0}H_{0,0})$.
\end{longtable}
On the other hand,
\begin{longtable}{lcl}
$D([L_{-\alpha,-i},H_{\alpha,i}])$&$=$&$[D(L_{-\alpha,-i}),H_{\alpha,i}]+[L_{-\alpha,-i},D(H_{\alpha,i})]$\\
$ $&$=$&$[f_{-\alpha}L_{-\alpha,-i}+g_{-\alpha}H_{-\alpha,-i},H_{\alpha,i}]+[L_{-\alpha,-i},\phi_{\alpha,i}L_{\alpha,i}+\psi_{\alpha,i}H_{\alpha,i}]$\\
$
$&$=$&$-\alpha f_{-\alpha}H_{0,0}-2\alpha\phi_{\alpha,i}L_{0,0}-\alpha\psi_{\alpha,i}H_{0,0}$.
\end{longtable}
Then we have $\alpha(2\phi_{\alpha,i}-\phi_{0,0})L_{0,0}+\alpha(f_{-\alpha}+\psi_{\alpha,i}-\psi_{0,0})H_{0,0}=0$, thus
$$\phi_{\alpha,i}=\frac{1}{2}\phi_{0,0},\,f_{-\alpha}+\psi_{\alpha,i}=\psi_{0,0}$$ for all $0\neq\alpha\in\Gamma.$
Since $H_{0,0}\in C(\mathcal{L}(\Gamma))$, then $D(H_{0,0})\in C(\mathcal{L}(\Gamma))$. Thus $\phi_{0,0}=0$. Furthermore, $\phi_{\alpha,i}=0$ for all $\alpha\in\Gamma$ and all $i\in\mathbb{Z}$, we have
\begin{longtable}{lcl}
$D(H_{\alpha,i})$&$=$&$\psi_{\alpha,i}H_{\alpha,i}$\\
$ $&$=$&$(-f_{-\alpha}+\psi_{0,0})H_{\alpha,i}$\\
$ $&$=$&$(f_\alpha+\psi_{0,0})H_{\alpha,i}$.
\end{longtable}
Let $\psi_{0,0}=b$. Then $D(H_{\alpha,i})=(f_\alpha+b)H_{\alpha,i}$.
\end{proof}
\end{lemma}

\begin{theorem}
$\mathrm{Der}\mathcal{L}(\Gamma)=\mathrm{ad}(\mathcal{L}(\Gamma))+(\mathcal{D}_{\mathrm{Hom}_{\mathbb{Z}}(\Gamma,\mathbb{F}[t,t^{-1}])}\oplus\mathcal{D}_{\mathbf{g}(\Gamma)}\oplus\mathcal{D}_{\mathbb{F}[t,t^{-1}]}\oplus
\mathcal{D}_{\mathbb{F}[t,t^{-1}]\frac{d}{dt}})$.

\begin{proof}
(\romannumeral1) First we show that $(\mathrm{Der}\mathcal{L}(\Gamma))_0=\mathcal{D}_{\mathrm{Hom}_{\mathbb{Z}}(\Gamma,\mathbb{F}[t,t^{-1}])}+\mathcal{D}_{\mathbf{g}(\Gamma)}+\mathcal{D}_{\mathbb{F}[t,t^{-1}]}+\mathcal{D}_
{\mathbb{F}[t,t^{-1}]\frac{d}{dt}}$.\\

In fact, by Lemma 3 and Lemma 4, on one hand, take any $D\in(\mathrm{Der}\mathcal{L}(\Gamma))_0$, for all $D_\phi\in\mathcal{D}_{\mathrm{Hom}_\mathbb{Z}(\Gamma,\mathbb{F}[t,t^{-1}])},\,D_g\in\mathcal{D}_{\mathbf{g}(\Gamma)},\,D_b\in\mathcal{D}_{\mathbb{F}[t,t^{-1}]},\,D^\rho\in\mathcal{D}_{\mathbb{F}[t,t^{-1}]\frac{d}{dt}}$, such that
\begin{longtable}{lcl}
$(D-D^\rho)(L_{\alpha,i})$&$=$&$f_\alpha L_{\alpha,i}+g_\alpha H_{\alpha,i}$,\\
$(D-D^\rho)(H_{\alpha,i})$&$=$&$(f_\alpha+b)H_{\alpha,i}$.
\end{longtable}
On the other hand, it follows  that
\begin{longtable}{lcl}
$(D_\phi+D_g+D_b)(L_{\alpha,i})$&$=$&$f_\alpha L_{\alpha,i}+g_\alpha H_{\alpha,i}$,\\
$(D_\phi+D_g+D_b)(H_{\alpha,i})$&$=$&$(f_\alpha+b)H_{\alpha,i}$.
\end{longtable}
Thus $D=D_\phi+D_g+D_b+D^\rho$.

(\romannumeral2) Next we show that $\mathcal{D}_{\mathrm{Hom}_{\mathbb{Z}}(\Gamma,\mathbb{F}[t,t^{-1}])}+\mathcal{D}_{\mathbf{g}(\Gamma)}+\mathcal{D}_{\mathbb{F}[t,t^{-1}]}+\mathcal{D}_{\mathbb{F}[t,t^{-1}]\frac{d}{dt}}$ is a direct sum.

In fact, let $$D=D_\phi+D_g+D_b+D^\rho,$$
where $\phi\in \mathrm{Hom}_{\mathbb{Z}}(\Gamma,\mathbb{F}[t,t^{-1}]),\,g\in\mathbf{g}(\Gamma),\,b\in\mathbb{F}[t,t^{-1}],\,\rho\in\mathbb{F}[t,t^{-1}]\frac{d}{dt}$. If $D=0$, to prove that $D_\phi=D_g=D_b=D^\rho=0$. In fact, we have
\begin{eqnarray}\nonumber
&&D(L_{\alpha,i})=0,\\ \nonumber
&&D(H_{\alpha,i})=0.
\end{eqnarray}
On the other hand,
\begin{eqnarray}\nonumber
&&(D_\phi+D_g+D_b+D^\rho)(L_{\alpha,i})=\phi(\alpha)L_{\alpha,i}+g_{\alpha}H_{\alpha,i}+L_\alpha\rho(t^i),\\ \nonumber
&&(D_\phi+D_g+D_b+D^\rho)(H_{\alpha,i})=\phi(\alpha)H_{\alpha,i}+bH_{\alpha,i}+H_{\alpha}\rho(t^i).
\end{eqnarray}
Then $D_\phi=D_g=D_b=D^\rho=0$.
\end{proof}
\end{theorem}

\section{$2$-local derivations on $\mathcal{L}(\Gamma)$}

Recall that a map $\Delta:\mathcal{L}\rightarrow\mathcal{L}$ (not linear in general) is called a $2$-local derivation if for every $x,y\in\mathcal{L}$, there exists a derivation $\Delta_{x,y}:\mathcal{L}\rightarrow\mathcal{L}$ such that $\Delta(x)=\Delta_{x,y}(x)$ and $\Delta(y)=\Delta_{x,y}(y)$ (see \cite{ref3}).

By Theorem 5, we have the following result.
\begin{corollary}
Let $\Delta$ be a $2$-local derivation on $\mathcal{L}(\Gamma)$. For all $x,y\in\mathcal{L}(\Gamma)$, there exist derivations $\Delta_{x,y}$ on $\mathcal{L}(\Gamma)$ such that $\Delta(x)=\Delta_{x,y}(x)\,,\,\Delta(y)=\Delta_{x,y}(y)$ and it can be written as $$\Delta_{x,y}=\mathrm{ad}(\sum_{\substack{\alpha\in\Gamma\\i\in\mathbb{Z}}}(a_{\alpha,i}(x,y)L_{\alpha,i}+b_{\alpha,i}(x,y)H_{\alpha,i}))+\lambda_1(x,y)D_\phi+\lambda_2(x,y)D_g+\lambda_3(x,y)D_b+\lambda_4(x,y)D^\rho.$$
where $a_{\alpha,i},b_{\alpha,i},\,\lambda_{i}$ for $i=1,2,3,4$ are complex-valued functions on $\mathcal{L}(\Gamma)\times\mathcal{L}(\Gamma)$, $D_\phi,D_g,D_b,D^\rho$ are given as previously.
\end{corollary}

\begin{lemma}
Let $\Delta$ be a $2$-local derivation on $\mathcal{L}(\Gamma)$. For any fixed $x\in\mathcal{L}(\Gamma)$, if $\Delta(L_{\alpha,i})=0$, $\alpha\in\Gamma,i\in\mathbb{Z}$, then
\begin{longtable}{lcl}
$\Delta_{L_{\alpha,i},x}$&$=$&$\mathrm{ad}(\sum_{\substack{\alpha\in\Gamma\\j\in\mathbb{Z}}}(a_{\alpha,j}(L_{\alpha,i},x)L_{\alpha,j}+b_{0,j}(L_{\alpha,i},x)H_{0,j}))+(1-\delta_{\alpha,0})\lambda_1(L_{\alpha,i},x)D_\phi$\\
$ $&&$+(1-\delta_{\alpha,0})\lambda_2(L_{\alpha,i},x)D_g+\lambda_3(L_{\alpha,i},x)D_b+(1-\delta_{i,0})\lambda_4(L_{\alpha,i},x)D^\rho$.
\end{longtable}

\begin{proof}
If $\Delta(L_{\alpha,i})=0$, then we have
\begin{longtable}{lcl}
$\Delta(L_{\alpha,i})$&$=$&$\Delta_{L_{\alpha,i},x}(L_{\alpha,i})$\\
$ $&$=$&$[\sum_{\substack{\beta\in\Gamma\\j\in\mathbb{Z}}}(a_{\beta,j}(L_{\alpha,i},x)L_{\beta,j}+b_{\beta,j}(L_{\alpha,i},x)H_{\beta,j}),L_{\alpha,i}]+\lambda_1(L_{\alpha,i},x)D_\phi(L_{\alpha,i})$\\ $ $&&$+\lambda_2(L_{\alpha,i},x)D_g(L_{\alpha,i})+\lambda_3(L_{\alpha,i},x)D_b(L_{\alpha,i})+\lambda_4(L_{\alpha,i},x)D^\rho(L_{\alpha,i})$\\ $ $&$=$&$\sum_{\substack{\beta\in\Gamma\\j\in\mathbb{Z}}}((\beta-\alpha)a_{\beta,j}(L_{\alpha,i},x)L_{\alpha+\beta,i+j}+\beta b_{\beta,j}(L_{\alpha,i},x)H_{\alpha+\beta,i+j})$\\
$ $&&$+\phi(\alpha)\lambda_1(L_{\alpha,i},x)L_{\alpha,i}+g_{\alpha}\lambda_2(L_{\alpha,i},x)L_{\alpha,i}+i\lambda_4(L_{\alpha,i},x)L_{\alpha,i-1}$\\
$ $&$=$&$0$.
\end{longtable}
Then $a_{\beta,j}(L_{\alpha,i},x)=0$ for $\alpha\neq\beta$, $b_{\beta,j}(L_{\alpha,i},x)=0$ for all $\beta\neq0$, $\lambda_1(L_{\alpha,i},x)=0$ for $\alpha\neq0$, $\lambda_2(L_{\alpha,i},x)=0$ for $\alpha\neq0$, $\lambda_4(L_{\alpha,i},x)=0$ for $i\neq0.$ Thus, the conclusion holds.
\end{proof}
\end{lemma}

\begin{lemma}
Let $\Delta$ be a $2$-local derivation on $\mathcal{L}(\Gamma)$ such that $\Delta(L_{0,i})=\Delta(L_{1,j})=0$. Then $\Delta(L_{\alpha,k})=0,$ for all $\alpha\in\Gamma$ and all $k\in\mathbb{Z}$.

\begin{proof}
Since $\Delta(L_{0,i})=\Delta(L_{1,j})=0$, then we can assume that
\begin{longtable}{lcl}
$\Delta_{L_{\gamma,z},x}$&$=$&$\mathrm{ad}(a_{\gamma,l}(L_{\gamma,z},x)L_{\gamma,l}+b_{0,l}(L_{\gamma,z},x)H_{0,l})+(1-\delta_{\gamma,0})\lambda_1(L_{\gamma,z},x)D_\phi$\\
$ $&&$+(1-\delta_{\gamma,0})\lambda_2(L_{\gamma,z},x)D_g+\lambda_3(L_{\gamma,z},x)D_b+(1-\delta_{z,0})\lambda_4(L_{\gamma,z},x)D^\rho$,
\end{longtable}
where $L_{\gamma,z}=\{L_{0,i},L_{1,j}\},\,x\in\mathcal{L}(\Gamma)$. Let $x=L_{\alpha,k}$.  Then we have
\begin{longtable}{lcl}
$\Delta(L_{\alpha,k})$&$=$&$\Delta_{L_{0,i},L_{\alpha,k}}(L_{\alpha,k})$\\
$ $&$=$&$[a_{0,l}(L_{0,i},L_{\alpha,k})L_{0,l}+b_{0,l}(L_{0,i},L_{\alpha,k})H_{0,l},L_{\alpha,k}]+(1-\delta_{0,0})\lambda_1(L_{0,i},L_{\alpha,k})D_\phi(L_{\alpha,k})$\\
$ $&&$+(1-\delta_{0,0})\lambda_2(L_{0,i},L_{\alpha,k})D_g(L_{\alpha,k})+\lambda_3(L_{0,i},L_{\alpha,k})D_b(L_{\alpha,k})$\\
$ $&&$+(1-\delta_{i,0})\lambda_4(L_{0,i},L_{\alpha,k})D^\rho(L_{\alpha,k})$\\
$ $&$=$&$-\alpha a_{0,l}(L_{0,i},L_{\alpha,k})L_{\alpha,l+k}+k(1-\delta_{i,0})\lambda_4(L_{0,i},L_{\alpha,k})L_{\alpha,k-1}$,\\
$\Delta(L_{\alpha,k})$&$=$&$\Delta_{L_{1,j},L_{\alpha,k}}(L_{\alpha,k})$\\
$ $&$=$&$[a_{1,l}(L_{1,j},L_{\alpha,k})L_{1,l}+b_{1,l}(L_{1,j},L_{\alpha,k})H_{0,l},L_{\alpha,k}]+(1-\delta_{1,0})\lambda_1(L_{1,j},L_{\alpha,k})D_\phi(L_{\alpha,k})$\\
$ $&&$+(1-\delta_{1,0})\lambda_2(L_{1,j},L_{\alpha,k})D_g(L_{\alpha,k})+\lambda_3(L_{1,j},L_{\alpha,k})D_b(L_{\alpha,k})$\\
$ $&&$+(1-\delta_{j,0})\lambda_4(L_{1,j},L_{\alpha,k})D^\rho(L_{\alpha,k})$\\
$ $&$=$&$(1-\alpha)a_{1,l}(L_{1,j},L_{\alpha,k})L_{\alpha+1,l+k}+\phi(\alpha)\lambda_1(L_{1,j},L_{\alpha,k})L_{\alpha,k}+g_\alpha\lambda_2(L_{1,j},L_{\alpha,k})L_{\alpha,k}$\\
$ $&&$+k(1-\delta_{j,0})\lambda_4(L_{1,j},L_{\alpha,k})L_{\alpha,k-1}$.
\end{longtable}
Comparing the formulas above, we have
\begin{longtable}{lcl}
$ $&&$a_{0,l}(L_{0,i},L_{\alpha,k})=a_{1,l}(L_{1,j},L_{\alpha,k})=\lambda_1(L_{1,j},L_{\alpha,k})=\lambda_2(L_{1,j},L_{\alpha,k})=0,$\\ $ $&&$\lambda_4(L_{0,i},L_{\alpha,k})=\lambda_4(L_{1,j},L_{\alpha,k})=0\,(k\neq0,\,i,j\neq0).$
\end{longtable}
When $k=0$ or $i,j=0$, we have $k(1-\delta_{i,0})=k(1-\delta_{j,0})=0$, thus $\Delta(L_{\alpha,k})=0$.
\end{proof}
\end{lemma}

\begin{lemma}
Let $\Delta$ be a $2$-local derivation on $\mathcal{L}(\Gamma)$ such that $\Delta(L_{\alpha,i})=0$, $\alpha\in\Gamma,i\in\mathbb{Z}$. Then $$\Delta(x)=0$$ for any $x=\sum_{\substack{\gamma\in\Gamma\\k\in\mathbb{Z}}}(m_{\gamma,k}L_{\gamma,k}+n_{\gamma,k}H_{\gamma,k})\in\mathcal{L}(\Gamma)$.

\begin{proof}
Given that $\Delta(L_{\alpha,i})=0$, by Lemma 7, we have
\begin{longtable}{lcl}
$\Delta(x)$&$=$&$\Delta_{L_{0,i},x}(x)$\\
$ $&$=$&$[a_{0,j}(L_{0,i},x)L_{0,j}+b_{0,j}(L_{0,i},x)H_{0,j},x]+\lambda_3(L_{0,i},x)D_b(x)+(1-\delta_{i,0})\lambda_4(L_{0,i},x)D^\rho(x)$\\
$ $&$=$&$\sum_{\substack{\gamma\in\Gamma k\in\mathbb{Z}}}(-\gamma m_{\gamma,k}a_{0,j}(L_{0,i},x)L_{\gamma,j+k}-\gamma n_{\gamma,k}a_{0,j}(L_{0,i},x)H_{\gamma,j+k})$\\
$ $&&$+n_{\gamma,k}b\lambda_3(L_{0,i},x)(H_{\gamma,k})+k(1-\delta_{i,0})\lambda_4(L_{0,i},x)(m_{\gamma,k}L_{\gamma,k-1}+n_{\gamma,k}H_{\gamma,k-1})$,\\
$\Delta(x)$&$=$&$\Delta_{L_{\alpha,i},x}(x)$\\
$ $&$=$&$[a_{\alpha,j}(L_{\alpha,i},x)L_{\alpha,j}+b_{0,j}(L_{\alpha,i},x)H_{0,j},x]+(1-\delta_{\alpha,0})\lambda_1(L_{\alpha,i},x)D_\phi(x)$\\
$ $&&$+(1-\delta_{\alpha,0})\lambda_2(L_{\alpha,i},x)D_g(x)+\lambda_3(L_{\alpha,i},x)D_b(x)+(1-\delta_{i,0})\lambda_4(L_{\alpha,i},x)D^\rho(x)$\\
$ $&$=$&$\sum_{\substack{\gamma\in\Gamma\\k\in\mathbb{Z}}}((\alpha-\gamma)m_{\gamma,k}a_{\alpha,j}(L_{\alpha,i},x)L_{\alpha+\gamma,j+k}-\gamma n_{\gamma,k}a_{\alpha,j}(L_{\alpha,i},x)H_{\alpha+\gamma,j+k})$\\
$ $&&$+\phi(\gamma)(1-\delta_{\alpha,0})\lambda_1(L_{\alpha,i},x)(m_{\gamma,k}L_{\gamma,k}+n_{\gamma,k}H_{\gamma,k})+g_\alpha(1-\delta_{\alpha,0})\lambda_2(L_{\alpha,i},x)(m_{\gamma,k}L_{\gamma,k})$\\
$ $&&$+b\lambda_3(L_{\alpha,i},x)(n_{\gamma,k}H_{\gamma,k})+k(1-\delta_{i,0})\lambda_4(L_{\alpha,i},x)(m_{\gamma,k}L_{\gamma,k-1}+n_{\gamma,k}H_{\gamma,k-1})$.
\end{longtable}
If there exists $m_{\alpha,i}\neq0$, taking sufficiently many $\alpha\in\Gamma,\,i\in\mathbb{Z}$, then we have $a_{\alpha,j}(L_{\alpha,i},x)=a_{0,j}(L_{\alpha,i},x)=0$ for all $i,j\in\mathbb{Z}$. Similarly $\lambda_3(L_ {0,i},x)=\lambda_3(L_{\alpha,i},x)=0$, and $\lambda_1(L_{\alpha,i},x)=\lambda_2(L_{\alpha,i},x)=0 \,(\alpha\neq0),\,\lambda_4(L_{\alpha,i},x)=0\,(i\neq0)$. When $\alpha=0$, we have $\phi(\alpha)(1-\delta_{\alpha,0})$ and $g_\alpha(1-\delta_{\alpha,0})=0$. When $i=0$, we have $i(1-\delta_{i,0})=0$.
\end{proof}
\end{lemma}

\begin{theorem}
Every $2$-local derivation on $\mathcal{L}(\Gamma)$ is a derivation.

\begin{proof}
Let $\Delta$ be a $2$-local derivation on $\mathcal{L}(\Gamma)$. Let $\Delta_{L_{0,i},L_{1,j}}$ be a derivation on $\mathcal{L}(\Gamma)$ satisfies
\begin{longtable}{lcl}
$\Delta(L_{0,i})$&$=$&$\Delta_{L_{0,i},L_{1,j}}(L_{0,i})$,\\
$\Delta(L_{1,j})$&$=$&$\Delta_{L_{0,i},L_{1,j}}(L_{1,j})$.
\end{longtable}
Let $\Delta(1)=\Delta-\Delta_{L_{0,i},L_{1,j}}$. Then $\Delta(1)$ is a $2$-local derivation satisfying $$\Delta(1)(L_{0,i})=\Delta(1)(L_{1,j})=0.$$ By Lemma 8, we have $\Delta(1)(L_{\alpha,k})=0$ and by Lemma 9, we have $\Delta(1)(x)=0,\,\mbox{for all}\, x\in\mathcal{L}(\Gamma)$. Then $\Delta(1)=\Delta-\Delta_{L_{0,i},L_{1,j}}\equiv0$. Furthermore, $\Delta=\Delta_{L_{0,i},L_{1,j}}$ is a derivation on $\mathcal{L}(\Gamma)$.
\end{proof}
\end{theorem}

\section{Biderivations on $\mathcal{L}(\Gamma)$}

\subsection{Biderivations on $\mathcal{L}(\Gamma)$}
Recall the definitions of biderivations and inner biderivations. Let $L$ be an algebra and $f$ be a bilinear map on $L\times L\rightarrow L$. If it
is a derivation with respect to both components, meaning that for all $x\in L$, there exists linear maps $\phi_x, \psi_x$ from $L$ to itself, and $\phi_x=f(x,\cdot), \psi_x=f(\cdot,x)$ are derivations on $L$, satisfying that
\begin{eqnarray}\nonumber
f(x\circ y,z)&=&x\circ f(y,z)+f(x,z)\circ y,\\ \nonumber
f(x,y\circ z)&=&f(x,y)\circ z+y\circ f(x,z).
\end{eqnarray}
Then $f$ is called a biderivation of $L$. A biderivation $f$ that satisfies $f(x,y)=\lambda[x,y],$ where $\lambda\in \mathbb{F}$, is called an inner biderivation.

It is easy to verify that the following lemma.

\begin{lemma}
$\mathcal{L}(\Gamma)$ is perfect.
\end{lemma}

\begin{lemma}(see \cite{ref10})
Let $L$ be perfect algebra and $f$ be a biderivation on $L$. If $\alpha\in C(L)$, then for all $x\in L$, $f(x,\alpha)=f(\alpha,x)=0$.
\end{lemma}

\begin{lemma}(see \cite{ref11})
Let $f$ be a biderivation on $L$. Then $$[f(x,y),[u,v]]= [[x,y],f(u,v)],\, x,y,u,v\in L.$$ In particular, $[f(x,y),[x,y]]=0,\, x,y\in L$.
\end{lemma}

\begin{lemma}(see \cite{ref11})
Let $f$ be a biderivation on $L$. If $[x,y]=0,\,x,y\in L$, then $f(x,y)\in C(L)$.
\end{lemma}

\begin{theorem}
Every biderivation of $\mathcal{L}(\Gamma)$ is inner.

\begin{proof}
Let $f$ be biderivation on $\mathcal{L}(\Gamma)$. We proceed by several steps.

(\romannumeral1) First, we show that there exists $\lambda\in\mathbb{F}$ such that $f(L_ {0,i},L_{m,j})\equiv\lambda[L_{0,i},L_{m,j}](\mathrm{mod}\mathbb{F}H_{0,i+j})$  for all $m\in\Gamma$ and $i,j\in\mathbb{Z}.$

If $m=0$, then $[L_{0,i},L_{0,j}]=0$, by Lemma 14 we conclude that $$f(L_{0,i},L_{0,j})\in C(\mathcal{L}(\Gamma))=\mathbb{F}H_{0,i+j}.$$

If $m\neq0$, let
$$f(L_{0,i},L_{m,j})=\sum_{\substack{n\in\Gamma\\i,j\in\mathbb{Z}}}(a_nL_{n,i+j}+b_nH_{n,i+j}).$$
By Lemma 13, we conclude
$$[[L_{0,i},L_{m,j}],f(L_{0,i},L_{m,j})]=0.$$
Since
$[-mL_{m,i+j},\sum_{\substack{n\in\Gamma\\i,j\in\mathbb{Z}}}(a_nL_{n,i+j}+b_nH_{n,i+j})]=0$ for all $n\in\Gamma$ and all $i,j\in\mathbb{Z},$ then it follows that
$$\sum_{\substack{n\in\Gamma}}(-ma_n(m-n))=0,\quad \sum_{\substack{n\in\Gamma}}mnb_n=0.$$
Furthermore, $a_n=0$ for $m\neq n$, $b_n=0$ for $n\neq0$, we conclude that
\begin{eqnarray}\nonumber
f(L_{0,i},L_{m,j})&\equiv&0(\mathrm{mod}\mathbb{F}H_{0,i+j}),\,n=0,\\ \nonumber
f(L_{0,i},L_{m,j})&\equiv&a_mL_{m,i+j},\,m=n,\\ \nonumber
f(L_{0,i},L_{m,j})&\equiv&0,\,m\neq n,n\neq0.
\end{eqnarray}
Thus $f(L_{0,i},L_{m,j})\equiv\lambda[L_{0,i},L_{m,j}](\mathrm{mod}\mathbb{F}H_{0,i+j})$ for all $m\in\Gamma$ and all $i,j\in\mathbb{Z}$.

(\romannumeral2) Secondly, we show that $$f(L_{0,i},H_{m,j})\equiv\lambda[L_{0,i},H_{m,j}](\mathrm{mod}\mathbb{F}H_{0,i+j})$$ for all $m\in\Gamma$ and $i,j\in\mathbb{Z}$.

If $m=0$, $H_{0,j}\in C(\mathcal{L}(\Gamma))$, by Lemma 12, $f(L_{0,i},H_{0,j})=0$.\\
When $m\neq0$, let $$f(L_{0,i},H_{m,j})=\sum_{\substack{n\in\Gamma\\i,j\in\mathbb{Z}}}(a_nL_{n,i+j}+b_nH_{n,i+j}).$$
Since
$$[[L_{0,j},H_{m,j}],f(L_{0,j},H_{m,j})]=0,$$
then we have
$[-mH_{m,i+j},\sum_{\substack{n\in\Gamma\\i,j\in\mathbb{Z}}}(a_nL_{n,i+j}+b_nH_{n,i+j})]=0$ for all $n\in\Gamma$ and all $i,j\in\mathbb{Z}.$
We conclude $\sum_{\substack{n\in\Gamma}}-m^2a_n=0$, thus $a_n=0$, then $$f(L_{0,i},H_{m,j})\equiv\sum_{\substack{n\in\Gamma\\i,j\in\mathbb{Z}}}b_nH_{n,i+j}.$$
By
$$[f(L_{0,i},H_{m,j}),[L_{0,i},L_{1,j}]]=[[L_{0,i},H_{m,j}],f(L_{0,i},L_{1,j})],$$
we have
$$[b_nH_{n,i+j},-L_{1,i+j}]=[-mH_{m,i+j},\lambda L_{1,i+j}],$$
which says
$$-nb_nH_{n+1,2(i+j)}=-m\lambda H_{m+1,2(i+j)}.$$
We have  $m=\sum_{\substack{n\in\Gamma}}n$, which gives  $f(L_{0,i},H_{m,j})\equiv\lambda[L_{0,i},H_{m,j}](\mathrm{mod}\mathbb{F}H_{0,i+j})$ for all $m\in\Gamma$ and all $i,j\in\mathbb{Z}$.

(\romannumeral3) Thirdly, we show that $f(L_{0,i},x)\equiv\lambda[L_{0,i},x](\mathrm{mod}\mathbb{F}H_{0,i+j})$ for all $x\in\mathcal{L}(\Gamma)$ and $i,j\in\mathbb{Z}$.

Let
$$x=\sum_{\substack{m\in\Gamma\\i,j\in\mathbb{Z}}}(a_mL_{m,j}+b_mH_{m,j}).$$ Then we have
\begin{eqnarray}\nonumber
f(L_{0,i},x)&=&\sum_{\substack{m\in\Gamma\\j\in\mathbb{Z}}}\left\{a_mf(L_{0,i},L_{m,j})+b_mf(L_{0,i},H_{m,j})\right\}\\ \nonumber
&=&\sum_{\substack{m\in\Gamma\\j\in\mathbb{Z}}}\left\{a_m\lambda[L_{0,i},L_{m,j}]+b_m\lambda[L_{0,i},H_{m,j}]\right\}(\mathrm{mod}\mathbb{F}H_{0,i+j})\\ \nonumber
&=&\lambda[L_{0,i},\sum_{\substack{m\in\Gamma\\j\in\mathbb{Z}}}(a_mL_{m,j}+b_mH_{m,j})](\mathrm{mod}\mathbb{F}H_{0,i+j})\\ \nonumber
&=&\lambda[L_{0,i},x](\mathrm{mod}\mathbb{F}H_{0,i+j}).
\end{eqnarray}
It follows that  $f(L_{0,i},x)\equiv\lambda[L_{0,i},x](\mathrm{mod}\mathbb{F}H_{0,i+j})$ for all $x\in\mathcal{L}(\Gamma)$ and all $i,j\in\mathbb{Z}$.

(\romannumeral4)Then, we show that $f(x,y)\equiv\lambda[x,y](\mathrm{mod}\mathbb{F}H_{0,i+j})$ for all $x,y\in\mathcal{L}(\Gamma)$ and all $i,j\in\mathbb{Z}$.

Following
$$[f(x,y),[L_{0,i},z]]=[[x,y],f(L_{0,i},z)]$$ for all $z\in\mathcal{L}(\Gamma),$
we have
$$[f(x,y)-\lambda[x,y],[L_{0,i},z]]=0$$ for all $z\in\mathcal{L}(\Gamma),$
it follows that   $f(x,y)-\lambda[x,y]\in C(\mathcal{L}(\Gamma))$,
thus $f(x,y)\equiv\lambda[x,y](\mathrm{mod}\mathbb{F}H_{0,i+j})$.

(\romannumeral5)Finally, we show that $f(x,y)\equiv\lambda[x,y]$ for all $x,y\in\mathcal{L}(\Gamma)$.

Assume that
$$f(x,y)=\lambda[x,y]+\alpha(x,y)H_{0,i+j},$$
where $\alpha$ is a bilinear map on $\mathcal{L}(\Gamma)\times\mathcal{L}(\Gamma)\rightarrow\mathbb{F}$.
We have
\begin{eqnarray}\nonumber
f([x,y],z)&=&\lambda[[x,y],z]+\alpha([x,y],z)H_{0,i+j},\\ \nonumber
f([x,y],z)&=&[x,\lambda[y,z]+\alpha(y,z)H_{0,i+j}]+[\lambda[x,z]+\alpha(x,z)H_{0,i+j},y]\\ \nonumber
&=&\lambda([x,[y,z]]+\lambda[[x,z],y]).
\end{eqnarray}
Then $$\alpha([x,y],z)=0$$ for all $x,y,z\in\mathcal{L}(\Gamma).$
We conclude $\alpha(x,y)=0$ for all $x,y\in\mathcal{L}(\Gamma)$. Finally, we obtain that $f(x,y)\equiv\lambda[x,y]$ for all $x,y\in\mathcal{L}(\Gamma)$. That is, $f$ is an inner derivation on $\mathcal{L}(\Gamma)$.
\end{proof}
\end{theorem}

\subsection{Applications}

 If a linear map $\phi$ on Lie algebra $L$, which subjects to $[\phi(x),x]=0$ for all $x\in L$, then $\phi$ called a linear commuting map on $L$. If $f(x,y)=[\phi(x),y]=[x,\phi(y)]$, then $f$ is a biderivation on $L$ (see \cite{ref11}).

\begin{theorem}
A linear map $\phi$ on $\mathcal{L}(\Gamma)$ is commuting if and only if there exist $\lambda\in\mathbb{F}$ and a linear map $\tau:\,\mathcal{L}(\Gamma)\rightarrow C(\mathcal{L}(\Gamma))$ such that $\phi(x)=\lambda x+\tau(x)$ for all $x\in\mathcal{L}(\Gamma)$.

\begin{proof}
Since
\begin{eqnarray}\nonumber
[\phi(x),x]&=&[\lambda x+\tau(x),x]\\ \nonumber
&=&\lambda[x,x]+[\tau(x),x]\\ \nonumber
&=&0,
\end{eqnarray}
then $\phi$ is a linear commuting map.

Conversely, let $f(x,y)=[\phi(x),y]$ for all $x,y\in\mathcal{L}(\Gamma)$. Then $f$\ is a biderivation on $\mathcal{L}(\Gamma)$. By Theorem 15, we have
$$f(x,y)=[\phi(x),y]=\lambda[x,y], \lambda\in \mathbb{F}.$$
Since $\phi(x)-\lambda x\in C(\mathcal{L}(\Gamma))$, then we have a linear map $\tau:\mathcal{L}(\Gamma)\rightarrow C(\mathcal{L}(\Gamma))$ such that $\phi(x)=\lambda x+\tau(x)$ for all $x\in\mathcal{L}(\Gamma)$.
\end{proof}
\end{theorem}

 Let $(L,[,])$ be a Lie algebra. A commutative $post$-Lie algebra structure on $L$ is a bilinear product $x\cdot y$ on $L$ satisfying the following identifies:
\begin{eqnarray}\nonumber
&&x\cdot y=y\cdot x,\\ \nonumber
&&[x,y]\cdot z=x\cdot(y\cdot z)-y\cdot(x\cdot z),\\ \nonumber
&&x\cdot[y,z]=[x\cdot y,z]+[y,x\cdot z],\, \mbox{for all}\, x,y,z\in L.
\end{eqnarray}
We also call $(\mathcal{L}(\Gamma),[,],\cdot)$ is a commutative $post$-Lie algebra (see \cite{ref12}).

\begin{lemma}(see \cite{ref12})
Let $(L,[,],\cdot)$ be a commutative $post$-Lie algebra. If $f:L\times L\rightarrow L$ subjects to $f(x,y)=x\cdot y,\, x,y\in L$, then $f$ is a biderivation on $L$.
\end{lemma}

\begin{theorem}
Any commutative $post$-Lie algebra structure on $\mathcal{L}(\Gamma)$ is trivial.

\begin{proof}
Assume $(\mathcal{L}(\Gamma),[,],\cdot)$ is a commutative $post$-Lie algebra. By Theorem 15 and Lemma 17, there exists $\lambda\in\mathbb{F}$ such that
$$f(x,y)=x\cdot y=\lambda[x,y].$$
Since $x\cdot y=y\cdot x$, we have
$$\lambda[x,y]=\lambda[y,x].$$
Thus $2\lambda[x,y]=0$, we conclude $\lambda=0$. That is $f(x,y)=x\cdot y=0$.
\end{proof}
\end{theorem}

\section{Automorphism groups on $\mathcal{L}(\Gamma)$}

Let $L=\mathrm{span}\{L_{\alpha ,i}\,|\,\alpha\in\Gamma,i\in\mathbb{Z}\},\,H=\mathrm{span}\{H_{\alpha ,i}\,|\,\alpha\in\Gamma,i\in\mathbb{Z}\}.$
\begin{lemma}
$H$ is a unique maximal ideal on $\mathcal{L}(\Gamma)$.

\begin{proof}
By the definition of $\mathcal{L}(\Gamma)$, $H$ is an ideal of $\mathcal{L}(\Gamma)$. Since $\mathcal{L}(\Gamma)/H$ is isomorphic to the generalized centerless loop Virasoro algebra. The centerless loop Virasoro algebra is a simple Lie algebra, then $H$ is a maximal ideal. Next, we prove $H=I$. Let $I$ be another maximal ideal of $\mathcal{L}(\Gamma)$. For all $0\neq x\in I$, let $x=l+h$, where $l\in L,\,h\in H$. If $l\neq0$, then $I=\mathcal{L}(\Gamma)$, which is a contradiction. Thus $l=0$.
\end{proof}
\end{lemma}

\begin{lemma}
Let $\mathcal{LV}(\Gamma)$ and $\mathcal{LV}(\Gamma^\prime)$ be two generalized centerless loop Virasoro algebra. Then $\mathcal{LV}(\Gamma)\cong\mathcal{LV}(\Gamma^\prime)$ if and only if there exists $a\in\mathbb{F}^\ast$ such that $a\Gamma^\prime=\Gamma$, where $\mathbb{F}^\ast=\mathbb{F}\setminus\{0\}$.

\begin{proof}
Notice that the ``only if'' part is obvious. We now prove the ``if'' part.

Given that $\mathcal{LV}(\Gamma)\cong\mathcal{LV}(\Gamma^\prime)$, let $\theta:\mathcal{LV}(\Gamma)\rightarrow\mathcal{LV}(\Gamma^\prime)$ be a Lie algebra isomorphism. Since $\mathbb{F}L_{0,0},\,\mathbb{F}L^\prime_{0,0}$ are the Cartan subalgebras of $\mathcal{LV}(\Gamma)$ and $\mathcal{LV}(\Gamma^\prime)$, respectively, then for all $a\in\mathbb{F}^\ast$ such that $\theta(L_{0,0})=aL^\prime_{0,0}.$
\begin{eqnarray}\nonumber
\theta([L_{0,0},L_{\alpha,j}])&=&\theta(-\alpha L_{\alpha,j})\\ \nonumber
&=&-\alpha\theta(L_{\alpha,j}),\\ \nonumber
\theta([L_{0,0},L_{\alpha,j}])&=&[\theta(L_{0,0}),\theta(L_{\alpha,j})]\\ \nonumber
&=&a[L^\prime_{0,0},\theta(L_{\alpha,j})].
\end{eqnarray}
Thus $[L^\prime_{0,0},\theta(L_{\alpha,j})]=-a^{-1}\alpha\theta(L_{\alpha,j})$, that is $-a^{-1}\alpha\in\Gamma^\prime$, we conclude $\Gamma\subset a\Gamma^\prime$. Similarly, it can be proved that $a\Gamma^\prime\subset\Gamma$.
\end{proof}
\end{lemma}

\begin{lemma}
$\mathcal{L}(\Gamma)\cong\mathcal{L}(\Gamma^\prime)$ if and only if there exists $a\in\mathbb{F}^\ast$ such that $a\Gamma^\prime=\Gamma$.

\begin{proof}
Let $\theta:\mathcal{L}(\Gamma)\rightarrow \mathcal{L}(\Gamma^\prime)$ be a Lie algebra isomorphism and $I,I^\prime$ be the maximal ideals of $\mathcal{L}(\Gamma)$ and $\mathcal{L}(\Gamma^\prime)$, respectively. By Lemma 19, we have $\theta(I)=I^\prime$, and by Lemma 20, we have $\mathcal{L}(\Gamma)/I\cong \mathcal{L}(\Gamma^\prime)/I^\prime$, then there exists $a\in\mathbb{F}^\ast$ such that $a\Gamma^\prime=\Gamma$.

Conversely, if $a\Gamma^\prime=\Gamma$, define a mapping by
\begin{eqnarray}\nonumber
\theta: \mathcal{L}(\Gamma)&\rightarrow&\mathcal{L}(\Gamma^\prime)\\ \nonumber
L_{\alpha,i}&\longmapsto&L_{\frac{\alpha}{a},i}\\ \nonumber
H_{\alpha,i}&\longmapsto&H_{\frac{\alpha}{a},i},\,\alpha\in\Gamma,\,i\in\mathbb{Z}.
\end{eqnarray}
Then $\theta$ is an isomorphism obviously.
\end{proof}
\end{lemma}

For all $\theta\in\mathrm{Aut}\mathcal{L}(\Gamma)$, considering $\theta$ act on the basis of $\mathcal{L}(\Gamma)$ for all $\alpha\in\Gamma,\,i\in\mathbb{Z}$, assume $f_{\alpha,i}\in\mathbb{F}[t,t^{-1}],$
\begin{eqnarray}
\theta(L_{\alpha,i})&=&L_{\frac{\alpha}{a}}f_{\alpha,i},\\
\theta(H_{\alpha,i})&=&H_{\frac{\alpha}{a}}f_{\alpha,i}.
\end{eqnarray}
Thus
\begin{eqnarray}\nonumber
\theta([L_{\alpha,i},L_{\beta,j}])&=&\theta((\alpha-\beta)L_{\alpha+\beta,i+j})\\ \nonumber
&=&(\alpha-\beta)L_{\frac{\alpha+\beta}{a}}f_{\alpha+\beta,i+j},\\ \nonumber
\theta([L_{\alpha,i},L_{\beta,j}])&=&[\theta(L_{\alpha,i}),\theta(L_{\beta,j})]\\ \nonumber
&=&[L_{\frac{\alpha}{a}}f_{\alpha,i},L_{\frac{\beta}{a}}f_{\beta,j}]\\ \nonumber
&=&\frac{\alpha-\beta}{a}L_{\frac{\alpha+\beta}{a}}f_{\alpha,i}f_{\beta,j},\\ \nonumber
\theta([L_{\alpha,i},H_{\beta,j}])&=&\theta(-\beta H_{\alpha+\beta,i+j})\\ \nonumber
&=&-\beta H_{\frac{\alpha+\beta}{a}}f_{\alpha+\beta,i+j},\\ \nonumber
\theta([L_{\alpha,i},H_{\beta,j}])&=&[\theta(L_{\alpha,i}),\theta(H_{\beta,j})]\\ \nonumber
&=&[L_{\frac{\alpha}{a}}f_{\alpha,i},H_{\frac{\beta}{a}}f_{\beta,j}]\\ \nonumber
&=&-\frac{\beta}{a}H_{\frac{\alpha+\beta}{a}}f_{\alpha,i}f_{\beta,j}.
\end{eqnarray}
By comparing the above expressions, we have:
\begin{eqnarray}\nonumber
&&a(\alpha-\beta)f_{\alpha+\beta,i+j}=(\alpha-\beta)f_{\alpha,i}f_{\beta,j},\\ \nonumber
&&a\beta f_{\alpha+\beta,i+j}=\beta f_{\alpha,i}f_{\beta,j},\,\mbox{for all}\,\alpha,\beta\in\Gamma\, \mbox{and all}\, i,j\in\mathbb{Z}.
\end{eqnarray}

\begin{lemma}
Let $a\in\mathbb{F}^\ast,\,f_{\alpha,i},f_{\beta,j}\in\mathbb{F}[t,t^{-1}]$. Then $af_{\alpha+\beta,i+j}=f_{\alpha,i}f_{\beta,j}$ for all $\alpha,\beta\in\Gamma$ and all $i,j\in\mathbb{Z}$.

\begin{proof}
The proof is similar with that in Lemma 3 and Lemma 4.
\end{proof}
\end{lemma}

Let $\mathbb{F}[t,t^{-1}]^\ast=\{at^i\,|\,a\in\mathbb{F}^\ast,i\in\mathbb{Z}\}$. By Lemma 22, we have $f_{0,0}$ is a non-zero element on $\mathbb{F}[t,t^{-1}]$, then $f_{0,0}\in\mathbb{F}[t,t^{-1}]^\ast$. Thus, we can write $f_{\alpha,i}$ in the following form:
\begin{eqnarray}\nonumber
&&f_{\alpha,i}:=a\mu(\alpha,i)t^{\epsilon(\alpha,i)},\,\mbox{for all}\,\alpha\in\Gamma \,\mbox{and all}\,i\in\mathbb{Z},
\end{eqnarray}
where $\mu(\alpha,i)\in\mathbb{F}^\ast,\,\epsilon(\alpha,i)\in\mathbb{Z}$, thus for all $\alpha_1,\alpha_2\in\Gamma$ and all $i,j\in\mathbb{Z}$,
\begin{eqnarray}\nonumber
&&\mu(\alpha_1,i)\mu(\alpha_2,j)=\mu(\alpha_1+\alpha_2,i+j),\\ \nonumber
&&\epsilon(\alpha_1,i)+\epsilon(\alpha_2,j)=\epsilon(\alpha_1+\alpha_2,i+j).
\end{eqnarray}
Then, (5.1) and (5.2) can be written as follows:
\begin{eqnarray}
\theta(L_{\alpha,i})&=&a\mu(\alpha,i)L_{{\frac{\alpha}{a}},\epsilon(\alpha,i)},\\
\theta(H_{\alpha,i})&=&a\mu(\alpha,i)H_{{\frac{\alpha}{a}},\epsilon(\alpha,i)}.
\end{eqnarray}
Based on the above formulas, we obtain the following automorphisms on $\mathcal{L}(\Gamma)$.

\begin{lemma}
(\romannumeral1) For all $a\in A=\{a\in\mathbb{F}^\ast\,|\,a\Gamma=\Gamma\}$, the map is given by
\begin{eqnarray}\nonumber
\theta_a:\mathcal{L}(\Gamma)&\longrightarrow&\mathcal{L}(\Gamma)\\ \nonumber
L_{\alpha,i}&\longmapsto&aL_{\frac{\alpha}{a},i}\\ \nonumber
H_{\alpha,i}&\longmapsto&aH_{\frac{\alpha}{a},i},\,\mbox{for all}\,\alpha\in\Gamma,\, \mbox{and all}\,i\in\mathbb{Z}
\end{eqnarray}
is an automorphism on $\mathcal{L}(\Gamma)$.

(\romannumeral2) For all $\phi\in\mathrm{Hom}(\Gamma,\mathbb{Z})$, the map is given by
\begin{eqnarray}\nonumber
\theta_\phi:\mathcal{L}(\Gamma)&\longrightarrow&\mathcal{L}(\Gamma)\\ \nonumber
L_{\alpha,i}&\longmapsto&L_{\alpha,i+\phi(\alpha)}\\ \nonumber
H_{\alpha,i}&\longmapsto&H_{\alpha,i+\phi(\alpha)},\,\mbox{for all}\,\alpha\in\Gamma,\, \mbox{and all}\,i\in\mathbb{Z}
\end{eqnarray}
is an automorphism on $\mathcal{L}(\Gamma)$.

(\romannumeral3) For all $\chi\in\chi(\Gamma)$, $\chi(\Gamma)$ is the set consisting of group homomorphisms from $\Gamma$ to $\mathbb{F}^\ast$, the map is given by
\begin{eqnarray}\nonumber
\theta_\chi:\mathcal{L}(\Gamma)&\longrightarrow&\mathcal{L}(\Gamma)\\ \nonumber
L_{\alpha,i}&\longmapsto&\chi(\alpha)L_{\alpha,i}\\ \nonumber
H_{\alpha,i}&\longmapsto&\chi(\alpha)H_{\alpha,i},\,\mbox{for all}\,\alpha\in\Gamma,\, \mbox{and all}\,i\in\mathbb{Z}
\end{eqnarray}
is an automorphism on $\mathcal{L}(\Gamma)$.

(\romannumeral4) For all $\psi\in\mathrm{Aut}\mathbb{Z}=\{id,-id\}$, the map is given by
\begin{eqnarray}\nonumber
\theta_\psi:\mathcal{L}(\Gamma)&\longrightarrow&\mathcal{L}(\Gamma)\\ \nonumber
L_{\alpha,i}&\longmapsto&L_{\alpha,\psi(i)}\\ \nonumber
H_{\alpha,i}&\longmapsto&H_{\alpha,\psi(i)},\,\mbox{for all}\,\alpha\in\Gamma,\, \mbox{and all}\,i\in\mathbb{Z}
\end{eqnarray}
is an automorphism on $\mathcal{L}(\Gamma)$.

(\romannumeral5) Given $b\in\mathbb{F}^\ast$, the map is given by
\begin{eqnarray}\nonumber
\theta_b:\mathcal{L}(\Gamma)&\longrightarrow&\mathcal{L}(\Gamma)\\ \nonumber
L_{\alpha,i}&\longmapsto&b^iL_{\alpha,i}\\ \nonumber
H_{\alpha,i}&\longmapsto&b^iH_{\alpha,i},\,\mbox{for all}\,\alpha\in\Gamma,\, \mbox{and all}\,i\in\mathbb{Z}
\end{eqnarray}
is an automorphism on $\mathcal{L}(\Gamma)$.

\begin{proof}
Obviously.
\end{proof}
\end{lemma}

Recall the \cite{ref6}, we have the set $(A\times\mathrm{Hom}(\Gamma,\mathbb{Z})\times\chi(\Gamma)\times\mathrm{Aut}\mathbb{Z}\times\mathbb{F^\ast})$ is a multiplicative group for which $(1,0,0,1,1)$ is the identity element, it satisfies $$(a_1,\phi_1,\chi_1,\psi_1,b_1)\cdot(a_2,\phi_2,\chi_2,\psi_2,b_2)=(a_1a_2,\phi_1\iota_{a^{-1}_2}+\psi_1\phi_2,b^{\phi_2}_1+\chi_2+\chi_1\iota_{a_2^{-1}},\psi_1\psi_2,b^{\psi_2}_1b_2),$$
where $\iota_a$ is an automorphism on $\Gamma$, satisfying $\iota_a(\alpha)=a\alpha$. If $b^\phi\in\chi(\Gamma)$, then $(b^\phi)(\alpha)=b^{\phi(\alpha)}$, where $\alpha\in\Gamma,\,b\in\mathbb{F}^\ast,\,\phi\in\mathrm{Hom}(\Gamma,\mathbb{Z}).$

\begin{theorem}
$\mathrm{Aut}\mathcal{L}(\Gamma)\cong(A\times Hom(\Gamma,\mathbb{Z})\times\chi(\Gamma)\times\mathrm{Aut}\mathbb{Z}\times\mathbb{F}^\ast)$.

\begin{proof}
Let $\Psi$ be a map from $A\times\mathrm{Hom}(\Gamma,\mathbb{Z})\times\chi(\Gamma)\times\mathrm{Aut}\mathbb{Z}\times\mathbb{F}^\ast$ to $\mathrm{Aut}\mathcal{L}(\Gamma)$,
\begin{eqnarray}\nonumber
\Psi:A\times\mathrm{Hom}(\Gamma,\mathbb{Z})\times\chi(\Gamma)\times\mathrm{Aut}\mathbb{Z}\times\mathbb{F}^\ast&\longrightarrow&\mathrm{Aut}\mathcal{L}(\Gamma)\\ \nonumber
(a,\phi,\chi,\psi,b)&\longmapsto&\theta_a\theta_\phi\theta_\chi\theta_\psi\theta_b
\end{eqnarray}
Now we prove $\Psi$ is an isomorphism. For all $$(a_1,\phi_1,\chi_1,\psi_1,b_1),(a_2,\phi_2,\chi_2,\psi_2,b_2)\in A\times\mathrm{Hom}(\Gamma,\mathbb{Z})\times\chi(\Gamma)\times\mathrm{Aut}\mathbb{Z}\times\mathbb{F}^\ast,$$
on the one hand,
\begin{eqnarray}\nonumber
&&\Psi(a_1,\phi_1,\chi_1,\psi_1,b_1)\Psi(a_2,\phi_2,\chi_2,\psi_2,b_2)(L_{\alpha,i})\\ \nonumber
&=&\Psi(a_1,\phi_1,\chi_1,\psi_1,b_1)(\theta_{a_2}\theta_{\phi_2}\theta_{\chi_2}\theta_{\psi_2}\theta_{b_2})(L_{\alpha,i})\\ \nonumber
&=&\Psi(a_1,\phi_1,\chi_1,\psi_1,b_1)(\theta_{a_2}\theta_{\phi_2}\theta_{\chi_2}\theta_{\psi_2})(b_2^iL_{\alpha,i})\\ \nonumber
&=&\Psi(a_1,\phi_1,\chi_1,\psi_1,b_1)(\theta_{a_2}\theta_{\phi_2}\theta_{\chi_2})(b^i_2L_{\alpha,\psi_2(i)})\\ \nonumber
&=&\Psi(a_1,\phi_1,\chi_1,\psi_1,b_1)(\theta_{a_2}\theta_{\phi_2})(b_2^i\chi_2(\alpha)L_{\alpha,\psi_2(i)})\\ \nonumber
&=&\Psi(a_1,\phi_1,\chi_1,\psi_1,b_1)(\theta_{a_2})(b_2^i\chi_2(\alpha)L_{\alpha,\psi_2(i)+\phi_2(\alpha)})\\ \nonumber
&=&\Psi(a_1,\phi_1,\chi_1,\psi_1,b_1)(a_2b_2^i\chi_2(\alpha)L_{\frac{\alpha}{a_2},\psi_2(i)+\phi_2(\alpha)})\\ \nonumber
&=&(\theta_{a_1}\theta_{\phi_1}\theta_{\chi_1}\theta_{\psi_1}\theta_{b_1})(a_2b_2^i\chi_2(\alpha)L_{\frac{\alpha}{a_2},\psi_2(i)+\phi_2(\alpha)})\\ \nonumber
&=&(\theta_{a_1}\theta_{\phi_1}\theta_{\chi_1}\theta_{\psi_1})(a_2b_1^{\psi_2(i)+\phi_2(\alpha)}b_2^i\chi_2(\alpha)L_{\frac{\alpha}{a_2},\psi_2(i)+\phi_2(\alpha)})\\ \nonumber
&=&(\theta_{a_1}\theta_{\phi_1}\theta_{\chi_1})(a_2b_1^{\psi_2(i)+\phi_2(\alpha)}b_2^i\chi_2(\alpha)L_{\frac{\alpha}{a_2},\psi_1(\psi_2(i)+\phi_2(\alpha))})\\ \nonumber
&=&(\theta_{a_1}\theta_{\phi_1})(a_2b_1^{\psi_2(i)+\phi_2(\alpha)}b_2^i\chi_1(\frac{\alpha}{a_2})\chi_2(\alpha)L_{\frac{\alpha}{a_2},\psi_1(\psi_2(i)+\phi_2(\alpha))})\\ \nonumber
&=&(\theta_{a_1})(a_2b_1^{\psi_2(i)+\phi_2(\alpha)}b_2^i\chi_1(\frac{\alpha}{a_2})\chi_2(\alpha)L_{\frac{\alpha}{a_2},\psi_1(\psi_2(i)+\phi_2(\alpha))+\phi_1(\frac{\alpha}{a_2})})\\ \nonumber
&=&a_1a_2b_1^{\psi_2(i)+\phi_2(\alpha)}b_2^i\chi_1(\frac{\alpha}{a_2})\chi_2(\alpha)L_{\frac{\alpha}{a_1a_2},\psi_1(\psi_2(i)+\phi_2(\alpha))+\phi_1(\frac{\alpha}{a_2})}.
\end{eqnarray}
On the other hand,
\begin{eqnarray}\nonumber
&&\Psi((a_1,\phi_1,\chi_1,\psi_1,b_1)\cdot(a_2,\phi_2,\chi_2,\psi_2,b_2))(L_{\alpha,i})\\ \nonumber
&=&\Psi(a_1a_2,\phi_1\iota_{a^{-1}_2}+\psi_1\phi_2,b^{\phi_2}_1+\chi_2+\chi_1\iota_{a_2^{-1}},\psi_1\psi_2,b^{\psi_2}_1b_2)(L_{\alpha,i})\\ \nonumber
&=&a_1a_2(b_1^{\psi_2}b_2)^ib_1^{\phi_2(\alpha)}\chi_2(\alpha)\chi(\frac{\alpha}{a_2})L_{\frac{\alpha}{a_1a_2},(\psi_1\psi_2)(i)+(\phi_1\iota_{a^{-1}_2}+\psi_1\phi_2)(\alpha)}\\ \nonumber
&=&a_1a_2b_1^{\psi_2(i)+\phi_2(\alpha)}b_2^i\chi_1(\frac{\alpha}{a_2})\chi_2(\alpha)L_{\frac{\alpha}{a_1a_2},\psi_1(\psi_2(i)+\phi_2(\alpha))+\phi_1(\frac{\alpha}{a_2})}.
\end{eqnarray}

Similarly, act on $H_{\alpha,i}$, we conclude that $$\Psi(a_1,\phi_1,\chi_1,\psi_1,b_1)\Psi(a_2,\phi_2,\chi_2,\psi_2,b_2)(H_{\alpha,i})=\Psi((a_1,\phi_1,\chi_1,\psi_1,b_1)\cdot(a_2,\phi_2,\chi_2,\psi_2,b_2))(H_{\alpha,i}),$$
thus $\Psi$ is a homomorphism. Next we prove $\Psi$ is a bijective map.

Let $\theta_{a_1}\theta_{\phi_1}\theta_{\chi_1}\theta_{\psi_1}\theta_{b_1}=\theta_{a_2}\theta_{\phi_2}\theta_{\chi_2}\theta_{\psi_2}\theta_{b_2}$. Since $\theta_a,\theta_\phi,\theta_\chi,\theta_\psi,\theta_b$ are automorphisms on $\mathcal{L}(\Gamma)$, then $\theta_{a_1}=\theta_{a_2},\theta_{\phi_1}=\theta_{\phi_2},\,\theta_{\chi_1}=\theta_{\chi_2},\,\theta_{\psi_1}=\theta_{\psi_2}, \theta_{b_1}=\theta_{b_2}$,furthermore, we conclude $(a_1,\phi_1,\chi_1,\psi_1,b_1)=(a_2,\phi_2,\chi_2,\psi_2,b_2)$, $\Psi$ is an injective map.

For all $\theta\in\mathrm{Aut}\mathcal{L}(\Gamma)$, by Lemma 21 and Lemma 22, we have $\theta(L_{\alpha,i})=a\mu(\alpha,i)L_{{\frac{\alpha}{a}},\epsilon(\alpha,i)}, \theta(H_{\alpha,i})=a\mu(\alpha,i)H_{{\frac{\alpha}{a}},\epsilon(\alpha,i)}$. Let $\mu(\alpha,0)=\chi(\alpha),\,\mu(0,i)=b(i),\,\epsilon(\alpha,0)=\phi(\alpha),\,\epsilon(0,i)=\psi(i)$. Then $\Psi(a,\phi,\chi,\psi,b)=\theta$, therefore $\Psi$ is a surjective map. In conclusion, $\Psi$ is an isomorphism.
\end{proof}
\end{theorem}

\end{document}